\documentclass[preprint,12pt]{elsarticle}
\usepackage{geometry}                
\usepackage{graphicx}
\usepackage{amssymb}
\usepackage{epstopdf}
\newtheorem{theorem}{Theorem}
\newtheorem{proposition}{Proposition}
\newtheorem{example}{Example}
\newtheorem{definition}{Definition}

\newtheorem{lemma}{Lemma}
\newtheorem{corollary}{Corollary}

\def\R{\Bbb R}

\def\spp{\vspace{5pt}  \noindent}

\newcommand{\eqref}[1]{{\rm (\ref{#1})}}

\input mssymb.tex

\usepackage{amssymb}

\journal{Differential Geometry and Applications}

\begin{document}

\begin{frontmatter}



\title{Minimum $L^\infty$ Accelerations in Riemannian Manifolds}


\author{Lyle Noakes}

\address{School of Mathematics and Statistics, The University of Western Australia, 35 Stirling Highway, Crawley, Western Australia 6009, Australia.\\
{\tt lyle@maths.uwa.edu.au}}

\begin{abstract}
Riemannian cubics are critical points for the $L^2$ norm of acceleration of curves in 
Riemannian manifolds $M$. In the present paper the $L^\infty$ norm replaces the $L^2$ norm, and a less direct argument is used to derive necessary conditions analogous to those in \cite{lylegreg}. The necessary conditions are examined when $M$ is a sphere or a bi-invariant Lie group.
\end{abstract}


\end{frontmatter}


\section{Introduction} Consider the task of moving a physical object between given configurations at given times, in a way that is {\em optimal} in some sense yet to be defined. 
\begin{example}\label{introex1} Let $t_0<t_1$ be given real numbers, and let $x_0,x_1\in E^3$ be given, where $E^m$ denotes Euclidean $m$-space. The task is to  move a point of unit mass from $x_0$ to $x_1$ along a curve $x:[t_0,t_1]\rightarrow E^3$ so as to minimise the total kinetic energy
$$K(x)~:=~\frac{1}{2}\int _{t_0}^{t_1}\Vert x^{(1)}(t)\Vert _E^2 ~dt$$
where $x^{(1)}$ denotes the first derivative of $x$ (the velocity), and $\Vert ~\Vert _E$ is the Euclidean norm. There is a unique minimising path, namely the uniform-speed interpolating line segment. 

\spp
Therefore, given real numbers $t_0<t_1<\ldots <t_n$ and points 
$x_0,x_1,\ldots ,x_n\in E^3$, the unique minimal energy interpolating curve is  a track-sum of $n$ line segments, traversed at uniform speed within each $[t_{j-1},t_j]$ where $j=1,2,\ldots ,n$. 

\spp
In practice there may be problems if an object is moved in this way, due to instantaneous changes of momentum at junctions $t_j$ for $j\not= 0,n$. So another measure of optimality might be used instead of $K$. A simple alternative is the mean squared norm of the acceleration 
$$J_2(x)~:=~\frac{1}{t_1-t_0}\int _{t_0}^{t_n}\Vert x^{(2)}(t)\Vert ^2~dt$$
where $x$ varies over piecewise-$C^2$ interpolating curves. As seen in \cite{deboorbook}, the {\em natural cubic spline} is the unique minimiser of $J_2$. 
\end{example}
For objects more complicated than a point mass, similar considerations may hold. 
\begin{example}\label{introex2} Suppose that we want to move a rigid body whose centre of mass is fixed through given configurations at times $t_0,t_1,\ldots ,t_n$. The trajectory becomes a curve $x:[t_0,t_n]\rightarrow M$, where 
$M$ is the $3$-dimensional rotation group $SO(3)$. The given configurations are points $x_j\in SO(3)$, and $K$ is now given by 
$$K(x)~=~ \frac{1}{2} \int _{t_0}^{t_1}\Vert x^{(1)}\Vert _{x(t)}^2~dt$$
where $\Vert ~\Vert _z$ is the Riemannian norm at $z\in M$. As in the Euclidean case, a  minimiser of $K$ is a track-sum of geodesic arcs traversed at uniform speed within each $[t_{j-1},t_j]$. 

\spp
Such a path might be impractical, for the same reasons as in Example \ref{introex1}. So, as before, we consider minimisers of 
$$J_2(x)~:=~\frac{1}{t_1-t_n}\int _{t_0}^{t_n}\Vert \nabla _tx^{(1)}\Vert _{x(t)}^2~dt$$
where $\nabla _tx^{(1)}$ is the {\em covariant acceleration} of the interpolant $x:[t_0,t_n]\rightarrow M$, and is adapted to the Riemannian structure of $M$. Minimisers of $J_2$, called  {\em Rieman cubic splines}, are $C^2$ track-sums of {\em Riemannian cubics} defined over the $t_{j-1},t_j$, with trivial covariant accelerations at $t_0$ and $t_n$. 

\spp
Riemannian cubics are $C^\infty$ curves that with more mathematical structure than either geodesics or cubic polynomials, and are given by the nonlinear ordinary differential equation \cite{lylegreg}
\begin{equation}
\nabla _t^3x^{(1)}+R(\nabla _tx^{(1)},x^{(1)})x^{(1)}~=~{\bf 0}
\end{equation}
where $R$ is the {\em Riemannian curvature} tensor field. 
\end{example}
Riemannian cubics are reviewed in  \cite{lyletomaszrev}.
More details, including recent results, can be found in \cite{jmp}, \cite{lylenonnull}, \cite{lyledual}, \cite{lyleqjm}, \cite{lyleSIAM}, \cite{pauley}, \cite{popiel}, \cite{vega}, \cite{crouch1}, \cite{crouch4}, \cite{giambo}. The present paper initiates an alternative line of enquiry, motivated by different engineering considerations. 

\spp
Whereas $J_2$ measures the mean squared norm of the force applied to the moving object, it may sometimes be more important to minimise the {\em maximum} norm of applied force, namely 
$$J_\infty (x)~:=~\max \{ \Vert \nabla _tx^{(1)}\Vert _{x(t)}:t\in [t_0,t_n]\}$$
where $x:[t_0,t_n]\rightarrow M$ varies over piecewise-$C^2$ interpolants. 
As seen in \S \ref{manyptssec}, in order to find minimising curves it suffices to take $n=1$ and impose additional constraints on velocities at endpoints. 

\spp
More precisely,  given $t_0<t_1$ together with $x_i\in M$ and $v_i\in TM_{x_i}$, a piecewise-$C^2$ curve $x:[t_0,t_1]\rightarrow M$ is {\em $F_{x_0,x_1,v_0,v_1}$-feasible} when 
$x(t_i)=x_i$ and $x^{(1)}(t_i)=v_i$ for $i=0,1$. A minimiser of $J_\infty$ on $F_{x_0,x_1,v_0,v_1}$-feasible curves is said to be {\em $F_{x_0,x_1,v_0,v_1}$-optimal}. 
Our main result, Theorem \ref{mainthm}, gives a necessary condition for a curve to be $F_{x_0,x_1,v_0,v_1}$-optimal. Examples are given in \S \ref{spheresec} and \S \ref{liegp} for the cases where  $M$ is a sphere and a bi-invariant Lie group respectively.  

\spp
\S \ref{2pts1velsec} introduces the related notion of {\em $F_{x_0,x_1,v_0}$-optimality}, for which an analogous necessary condition is given by  Theorem \ref{subthm1}. Then  Corollary \ref{subthm2} of \S \ref{manyptssec} says that, if $x$ minimises $J_\infty$ over piecewise-$C^2$ interpolants of many points, then a least one track-summand of $x$ satisfies the conclusions of Theorem \ref{mainthm} or Theorem \ref{subthm1}. 
\section{Optimal Curves}
Let $M$ be a $C^\infty$ manifold of dimension $m\geq 1$, with a $C^{\infty}$ Riemannian metric $\langle ~,~\rangle $. Let real numbers $t_0<t_1$be given. The derivative $x^{(1)}$ of a piecewise-$C^2$ curve 
$x:[t_0,t_1]\rightarrow M$ is piecewise-$C^1$, and the Levi-Civita covariant acceleration $\nabla _tx^{(1)}$ is a piecewise-continuous field defined along $x$. 
\begin{definition}  
A piecewise-$C^2$ curve $x:[t_0,t_1]\rightarrow M$ is  said to be {\em feasible} when 
$x(t_i)=x_i$ and $x^{(1)}(t_i)=v_i$ for $i=0,1$.  
\end{definition}
Let $F=F_{x_0,x_1,v_0,v_1}=F_{x_0,x_1,v_0,v_1,t_0,t_1}$ be the set of all feasible curves, and  define 
$J:F\rightarrow \R $
by $J(x):=\max \{\Vert \nabla _tx^{(1)}\Vert _{x(t)} :t\in [t_0,t_1]\}$ where, for $w\in M$,  $\Vert ~\Vert _w$ is the Riemannian norm on $TM_w$ associated with the Riemannian inner product $\langle ~,~\rangle _w$ . A minimizer of $J$ is said to be {\em optimal}. 
Our main result is a necessary condition for optimality, in terms of  
$$L(X)~:=~\nabla _t^2X+R(X,x^{(1)})x^{(1)}$$
where $R$ denotes Riemannian curvature, and $X$  a piecewise-$C^2$ field defined along $x$.
\begin{theorem}\label{mainthm} Let $x\in F$ be optimal. Then    
$\Vert  \nabla _tx^{(1)} \Vert _{x(t)}$ is constant. For some piecewise-$C^2$ function 
$\Phi :[t_0,t_1]\rightarrow [0,\infty )$ with $\Phi ^{-1}(0)$ discrete, we have 
\begin{equation}\label{eleq}L(\phi \nabla _tx^{(1)})~=~{\bf 0}\end{equation}
where $\phi (t):=\sqrt{\Phi (t)}$. Moreover $x$ and $\phi $ are $C^\infty$ except possibly where $\Phi (t)=0$. 
\end{theorem}
\begin{example}\label{eucex} Let $x$ be an optimal curve in an open convex subset $M$ of Euclidean $m$-space $E^m$. Then (\ref{eleq}) reads
$$\frac{d^2}{dt^2}(\phi x^{(2)})~=~{\bf 0}~\Longrightarrow ~\phi (t)x^{(2)}(t)~=~A+Bt$$ 
where $A,B\in E^m$. We have  
$z\phi (t)=\Vert A+Bt\Vert $ where $z$ is the constant length of $x^{(2)}$. 
\begin{itemize}
\item If $z=0$ then $A=B={\bf 0}$ and $x$ is an affine line segment.
\item If $A+Bt_2={\bf 0}$ for some $t_2\in [t_0,t_1]$ where $B\not= {\bf 0}$ then 
$$\phi (t)~=~\frac{\vert t-t_2\vert\Vert B\Vert  }{a}~\Longrightarrow ~x^{(2)}(t)~=~z~{\rm sign }(t-t_2)\frac{B}{\Vert B\Vert }~\hbox{ for }t\not= t_2.$$
So $x$ is a $C^1$ quadratic polynomial spline which is $C^\infty$ except possibly at $t_2$.  
\item If $A$ and $B$ are linearly independent then $\phi $ is never zero. By shifting $[t_0,t_1]$ we can suppose $\langle A,B\rangle =0$. Setting 
$\alpha :=\Vert A\Vert $ and $\beta :=\Vert B\Vert $, 
$$x^{(2)}(t)~=~z\frac{A+Bt}{\sqrt{\alpha ^2+\beta ^2t^2}}~\Longrightarrow $$
\begin{eqnarray*}
x(t)&=&z\frac{\beta t\log (\beta t+\sqrt{\alpha ^2+\beta ^2t^2})-\sqrt{\alpha  ^2+\beta ^2t^2}}
{\beta ^2}A+\\
&~&z\frac{\alpha ^2\log (\beta t+\sqrt{\alpha ^2+\beta ^2t^2})+\beta t\sqrt{\alpha ^2+\beta ^2t^2}}{2\beta ^3}B+Ct+D
\end{eqnarray*}
where $C,D\in E^m$ are constants. 
\end{itemize}
\end{example}
The conclusions of Example \ref{eucex} also follow from the argument preceding the proof of Theorem 4.1  in \S 4 of \cite{charlie}.  However \cite{charlie}, uses linear methods from classical approximation theory, building on \cite{deboor}, \cite{jerome}, \cite{glaeser} among others. These methods are inapplicable in a more general Riemannian setting. For instance, in \S \ref{spheresec} where $M$ is a sphere, or in \S \ref{liegp} where $M$ is a bi-invariant semisimple Lie group, there seems to be no alternative to Theorem \ref{mainthm}.

\spp
The proof of Theorem \ref{mainthm} is less direct than the derivations of the analogous Euler-Lagrange equations for geodesics \cite{hawking}, elastic curves \cite{singer},  and Riemannian cubics \cite{lylegreg}. As with these other variational problems, we proceed by considering {\em extremals}, where an optimal curve is necessarily extremal, but we need two distinct notions of extremal. 

\spp
In \S \ref{Linfsec}, so-called {\em $L^\infty$-extremals} are used to prove that $\Vert \nabla _tx^{(1)}(t)\Vert _{x(t)}$ is constant, and that optimal curves are {\em locally optimal}. After some preliminaries in \S \ref{localsec}, local optimality is used in \S \ref{reducesec} to reduce Theorem \ref{mainthm} to the special case where $M$ is an open subset of $\R ^m$ with a nonstandard Riemannian metric. For this special case, constancy of $\Vert \nabla _tx^{(1)}(t)\Vert _{x(t)}$ is used to identify $x$ with an extremal 
in the sense of {\em optimal control} for a control system with state variables in $\R ^{2m+1}$ and controls in the Euclidean unit sphere $S^{m-1}$. Then in \S \ref{localcalcsec} a local calculation based on the Pontryagin Maximum Principle \cite{ponty} is used to prove Theorem \ref{mainthm} for the special case, and this completes the proof of Theorem \ref{mainthm} in general. 

\spp
A more direct proof of Theorem \ref{mainthm} is possible for optimal curves 
whose image is contained in a single coordinate chart, based on the Pontryagin 
Principle for 
a slightly different optimal control  problem where $\Vert \nabla _tx^{(1)}\Vert _{x(t)}$ is not assumed to be constant. Then \S \ref{Linfsec} could be omitted, at the expense of generality. Another alternative might be to use the Pontryagin Principle for multiprocesses \cite{clarke}, but \S \ref{Linfsec} seems less complicated. 

\spp
It might be thought that an application of the geometric formulation of the Pontryagin Principle could replace the local calculations of \S \ref{localcalcsec}, and perhaps do away with \S \ref{Linfsec}. For this to work, the space of controls would need to vary from point to point, which is not envisaged in \cite{bloch}, \cite{barbero}, \cite{chang}. Possibly this might be done in the framework of \cite{agrachev}, but our local calculation has the advantage of being more elementary. 

\spp
Having proved Theorem \ref{mainthm}, we go on to investigate some consequences. In  
\S \ref{spheresec},  $M$ is taken to be the Euclidean unit $m$-sphere, and the differential equation (\ref{eleq}) 
is rewritten as the system  (\ref{xeq}), (\ref{Xeq}), without covariant derivatives. 
Figure \ref{xfig} in Example \ref{sphereex} shows a numerical solution for $x$. 

\spp
In \S \ref{liegp}, $M$ is taken to be  a Lie group $G$ with a bi-invariant Riemannian metric. The differential equation (\ref{eleq}) reduces and integrates to a differential equation (\ref{intreducedeleq}) in the Lie algebra ${\cal G}$,  in terms of a {\em Lie reductions} $V,X_L$. 
Another two conserved scalar quantities in terms of $V,X_L$ are also noted, and  
it is shown that generically $x$ can be recovered by quadrature from $V,X_L$.  
A class of so-called { null optimal curves} turn out to be the same as null Riemannian cubics \cite{jmp}, \cite{lyledual}, \cite{lyleSIAM} and Chapter 4 of \cite{pauley}. The curious geometry of reductions of non-null optimal curves in $SO(3)$ is illustrated by Figures \ref{vfig}, \ref{phfig}, \ref{vfig2} of Example \ref{so3ex}.  

\spp
In \S \ref{2pts1velsec}, $M$ is once more a general Riemannian manifold, and the notion of feasibility is relaxed, so that $x^{(1)}$ is specified only at one endpoint. Simple modifications to the proof of Theorem \ref{mainthm} result in an additional necessary condition, given in Theorem \ref{subthm1}. 

\spp
In \S \ref{manyptssec} Theorems \ref{mainthm}, \ref{subthm1} are adapted to a $J_\infty$ analogue of the situation in Example \ref{introex2}, where the $x(t_i)$ are specified at various $t_i$ and $x^{(1)}$ is unconstrained. 
As in Corollary \ref{subthm2}, a necessary condition for such a feasible curve 
$x$ to be optimal is that equation (\ref{eleq}) should hold along at least one arc, with an additional conditions when the arc is terminal. 

\section{$L^\infty$-Extremals}\label{Linfsec}
A piecewise-$C^2$ field $W:[t_0,t_1]\rightarrow TM$ along $x\in F$  
is {\em variational} when $W(t_i)=\nabla _tW\vert _{t=t_i}={\bf 0}$ for $i=0,1$. For $W$ variational, a {\em variation} $h\mapsto x_h\in F$ is given by  
$x_h(t):={\rm exp}_{x(t)}(hW(t))$, and
$$(\frac{\partial}{\partial h}x_h(t))\vert _{h=0}~=~W(t).$$
For $x\in F$ set $S_x:=\{ s\in [t_0,t_1]:\Vert \nabla _tx^{(1)}(t)\Vert _{s=t}=J(x)\} $. Then $S_x$ is nonempty and closed. 
\begin{definition} $x\in F$ is an {\em $L^\infty$-extremal} when, for any variational field $W$ along $x$, we have $\langle L(W),\nabla _tx^{(1)}\rangle  _{t=s_W}\geq 0$
for some $s_W\in S_x$.
\end{definition}  
\begin{lemma}\label{lem1} If $x$ is optimal then $x$ is an $L^\infty$-extremal. 
\end{lemma}

\spp
{\bf Proof:} From the definitions of $x_h$ and $L$, we find, for any variational field $W$,   
$$\frac{\partial }{\partial h}(\Vert \nabla _tx_h^{(1)}\Vert ^2)_{h=0}~=~
2\langle L(W),\nabla _t x^{(1)}\rangle .$$ 
Suppose $x\in F$ is not an $L^\infty$-extremal. Then, for some variational field $W$,  $\langle L(W),\nabla _t x^{(1)}\rangle  < 0$ at all $s\in S_x$. 
Therefore, and because $S_x$ is compact, for $h$ sufficiently small  we have $\Vert \nabla _t x_h^{(1)}\Vert ^2 _{t=s}<J(x)^2$ for all $s\in S_x$. Then $\Vert \nabla _t x_h^{(1)}\Vert ^2<J (x)^2$ at all $t$ in an open subset $U$ of $[t_0,t_1]$ containing $S_x$. In particular, if $U=[t_0,t_1]$ then $x$ is not optimal.

\spp
If $U\subset [t_0,t_1]$ then, because the nonempty set 
$[t_0,t_1]-U$ is compact, 
$\sup \{ \Vert \nabla _t x_h^{(1)}\Vert ^2:t\notin U\} $ is $\Vert \nabla _t x_h^{(1)}\Vert ^2 _{t=t_*}$ for some $t_*\notin U$. Because $S_x\subset U$, $\Vert \nabla _tx_h^{(1)}\Vert ^2 _{t=t_*}<J(x)^2$. So $\Vert \nabla _t x_h^{(1)}\Vert ^2<J (x)^2$ on $[t_0,t_1]$, and $x$ is not optimal. This proves Lemma \ref{lem1}. 

\begin{lemma}\label{lem2} If $x$ is an $L^\infty$-extremal then $\Vert \nabla _tx^{(1)}\Vert _{x(t)} $ is constant. 
\end{lemma}

\spp
{\bf Proof:} If $x$ is a geodesic the lemma holds trivially. For $x$ not a geodesic we have $J(x)>0$.  Arguing  by contradiction, suppose that the open set $U:=[t_0,t_1]-S_x$ is nonempty. Then $U$ contains an open interval of the form 
$(t_*-\delta ,t_*+\delta )$ where $t_*\in [t_0,t_1]$ and $\delta >0$. Define fields 
$\hat W_L:[t_0,t_*-\delta ]\rightarrow TM$ along $x\vert [t_0,t_*-\delta ]$ and 
$\hat W_R:[t_*+\delta ,t_1]\rightarrow TM$ along $x\vert [t_*+\delta ,t_1]$ 
by solving the linear initial and terminal value problems 
\begin{eqnarray*}L(\hat W_L)~=~-\nabla _tx^{(1)}\quad &\hbox{subject to }&\hat W_L(t_0)=\nabla _t\hat W_L\vert _{t=t_0}={\bf 0}\\ 
L(\hat W_R)~=~-\nabla _tx^{(1)}\quad &\hbox{subject to }&\hat W_R(t_1)=\nabla _t\hat W_R\vert _{t=t_1}={\bf 0}.\end{eqnarray*} 
Then $\hat W_L$ and $\hat W_R$ are the restrictions of a single variational field $W:[t_0,t_1]\rightarrow TM$ along $x$, for which $\langle L(W),\nabla _tx^{(1)}\rangle _{t=s}=-J (x)^2$ for all $s\in S_x$. Since $x$ is an $L^\infty$-extremal we have a contradiction. This proves Lemma \ref{lem2}. 

\spp
Incidentally, there is also a kind of converse. 
\begin{proposition}\label{prop1} Let $x\in F$ satisfy the necessary condition of Theorem \ref{mainthm}.  
Then $x$ is an $L^\infty$-extremal. 
\end{proposition}

\spp
{\bf Proof:} Set $X(t):=\phi (t)\nabla _tx^{(1)}(t)$. If $x$ is not an $L^\infty$-extremal then, for some variational field $W$, 
$$\langle \nabla _t^2W+R(W,x^{(1)})x^{(1)},X\rangle ~>~0\quad \hbox{on~}[t_0,t_1).$$
Using symmetries of the curvature tensor \cite{hawking} (2.27), and equation (\ref{eleq}), the left hand side is 
$$\langle \nabla _t^2W,X\rangle +\langle W,R(X,x^{(1)})x^{(1)}\rangle =
\langle \nabla _t ^2W,X\rangle -\langle W,\nabla _t^2X\rangle ~=~r^{(1)}(t)$$
where $r:[t_0,t_1]\rightarrow \R $ is given by $r(t):=\langle \nabla _tW,X\rangle -\langle W,\nabla _tX\rangle $.  So $r$ is strictly increasing. 

\spp
Because $W(t_i)=\nabla _tW\vert _{t=t_i}={\bf 0}$, $r(t_i)=0$ for $i=0,1$.   This is a contradiction because $r$ is strictly increasing. So $x$ is an $L^\infty$-extremal after all, and this proves Proposition \ref{prop1}.

\spp
Proposition \ref{prop1} is not used for the proof of Theorem \ref{mainthm}. Instead we need Lemma \ref{locallem} below, which asserts that an optimal curve is also {\em locally optimal}. 
This would be trivial if $J(x)$ was the integral of a non-negative function of $x$ and its derivatives, as with geodesics, elastic curves and Riemannian cubics. In the present situation a proof is needed. 
\begin{lemma}\label{locallem} Let $x\in F_{x_0,x_1,v_0,v_1,t_0,t_1}$ be optimal, and choose $s_0<s_1$ in $[t_0,t_1]$. For $i=0,1$ set $y_i=x(s_i)$ and $w_i=x^{(1)}(s_i)$. Then the restriction 
$y=x\vert [s_0,s_1]$ of $x$ to $[s_0,s_1]$ is optimal, considered as a curve in $F_{y_0,y_1,w_0,w_1,s_0,s_1}$.  
\end{lemma} 

\spp
{\bf Proof:} By Lemma \ref{lem2}, $J(y)=J(x)$. If $y$ is not optimal let $\bar y\in F_{y_0,y_1,w_0,w_1,s_0,s_1}$ be a curve with 
$\Vert \nabla _t\bar y^{(1)}\Vert _s<J(y)$ for all $s\in [s_0,s_1]$.  The track-sum $\bar x:[t_0,t_1]\rightarrow M$ of $x\vert [t_0,s_0]$, $\bar y$ and $x\vert [s_1,t_1]$ is piecewise-$C^2$. Indeed $\bar x\in F_{x_0,x_1,v_0,v_1,t_0,t_1}$. 

\spp
If $s_0=t_0$ and $s_1=t_1$ the Lemma holds trivially. So we can suppose without loss of generality that $t_2\notin [s_0,s_1]$ for some $t_2\in [t_0,t_1]$. So $J(\bar x)=\Vert \nabla _tx^{(1)}\vert _{t=t_2}\Vert _{x(t_2)}=J(x)$. So $\bar x$ is optimal. Since $\Vert \nabla _t\bar y^{(1)}\Vert _s<J(y)$ for all $s\in [s_0,s_1]$, 
$\Vert \nabla _t\bar x\Vert _{\bar x(t)}$ is nonconstant, and this contradicts Lemma \ref{lem2}. So $y$ is optimal after all, and Lemma \ref{locallem} is proved. 

\spp
Given $s_0=t_0<s_1<\ldots <s_j<\ldots s_n=t_1$, any curve $x\in F$ is a track-sum of its restrictions to the subintervals $[s_{j-1},s_j]$ . If $x$ is optimal then, by Lemma \ref{locallem}, so are the $x\vert [s_{j-1},s_j]$ with respect to the values and derivatives of $x$ at $s_{j-1}$ and $s_j$. So it suffices to prove Theorem \ref{mainthm} when $x$ maps into a coordinate chart of $M$. In this case, $M$ can be replaced 
by an open subset of $\R ^m$ equipped with some Riemannian 
metric $\langle ~,~\rangle $. From now on suppose this has been done.
\section{Local Geometry}\label{localsec}
Mainly to be clear in future on matters of notation, we briefly review some coordinate-based differential geometry, where $M$ is taken as an open subset of $\R ^m$. Readers may prefer to skip  
to Lemma \ref{lem3}, then \S \ref{reducesec} and onwards, referring back as necessary.  

\spp
For any $v\in M\subseteq \R ^m$ we have a possibly non-Euclidean inner product $\langle ~,~\rangle _v$ on $\R ^m=TM_v$. The associated norm 
is denoted by $\Vert ~\Vert _v$, and the { dual} $\tilde \omega \in \R ^m$ of a linear form $\omega \in (\R ^m)^*$ with respect to $\langle ~,~\rangle _v$ is given by 
$\omega (w) =\langle \tilde \omega ,w\rangle _v$ for any $w\in \R ^m$. 

\spp
Any vector $w\in \R ^m$ can be written in the form $w^ie_i$, where summation over $i$ is understood, where the $w^i\in \R$, and where $e_1,e_2,\ldots ,e_m$ are the standard basis elements of $\R ^m$. For any vector field $W$ on $M$, $W=W^ie_i$ where the $W^i:M\rightarrow \R$. The Levi-Civita covariant derivative $\nabla _YW$ of $W$ at $Y\in TM_v$ is the vector in $TM_v=\R ^m$ given by 
$$\nabla _YW~=~\frac{\partial W}{\partial x_i}Y^i+\Gamma _v(Y,W)$$
where, for each $v\in M$, $\Gamma _v:TM_v\times TM_v\rightarrow TM_v$ is the symmetric bilinear form given by 
$$\Gamma _v(Y,W)~:=~{\Gamma ^k_{ij}}(v)Y^iW^je_k$$
and the {\em Christoffel symbols} ${\Gamma ^k_{ij}}$ are $C^\infty$ functions for all $1\leq i,j,k\leq m$. The Riemannian metric determines the Christoffel symbols according to formula (2.26) of \cite{hawking}, and thereby the {\em Riemannian curvature} $R_v(X,Y)Z$ which is trilinear in $X,Y,Z\in TM_v$ and given by 
\begin{equation}\label{curveq}R(X,Y)Z~=~(\frac{\partial \Gamma ^i_{pj}}{\partial x_k}-\frac{\partial \Gamma ^i_{kj}}{\partial x_p}+\Gamma ^i_{kq}\Gamma ^q_{pj} - \Gamma ^i_{pq}\Gamma ^q_{kj} )X^kY^pZ^je_i\end{equation}
according to formula (2.20) of \cite{hawking}, where summation over $i,j,k,p,q$ is understood. 

\spp
Now, given $x\in F=F_{x_0,x_1,v_0,v_1}$, the velocity $y:=x^{(1)}$ is a vector field along $x$, and the {\em covariant accleration} of $x$ is given by the formula
$$\nabla _tx^{(1)}~:=~\nabla _{x^{(1)}(t)}y~=~y^{(1)}(t)+\Gamma (x(t))(y(t),y(t))$$
and an optimal curve is one for which the maximum Riemannian norm of the right hand side is minimised. By Lemmas \ref{lem1}, \ref{lem2}, it suffices to optimise over curves $x\in F$ for which $z:=\Vert \nabla _tx^{(1)}\Vert _{x(t)} $ is constant. Since Theorem \ref{mainthm} holds trivially when $z=0$, take $z>0$.

\spp
Let $S:=S^{m-1}\subset E^m$ be the unit sphere with respect to the {\em Euclidean} norm $\Vert ~\Vert _E$ and define $f:M\times S\rightarrow \R ^m$ by 
$$f(w,u)~:=~\frac{u~~}{\Vert u\Vert _w}.$$
\begin{lemma}\label{lem3} For any $u\in S$ and any $j=1,2,\ldots m$, 
$$\frac{\partial }{\partial w_j}f(w,u)=-\langle \Gamma (e_j,f(w,u)),f(w,u)\rangle _w f(w,u).$$
\end{lemma}

\spp
{\bf Proof of Lemma \ref{lem3}:} Because the Levi-Civita covariant derivative is compatible with the Riemannian metric, 
$$\Vert u\Vert _w\frac{\partial}{\partial w_j}\Vert u\Vert _w~=~\frac{1}{2}\frac{\partial }{\partial w_j}\langle u,u\rangle _w~=~\langle  \Gamma (e_j,u),u\rangle _w~\Longrightarrow ~\frac{\partial}{\partial w_j}(\Vert u\Vert _w)~=~\langle  \Gamma (e_j,f(w,u)),u\rangle _w$$ 
which proves Lemma \ref{lem3}.
\section{Local Reduction to Optimal Control}\label{reducesec}
For $M$ open in $\R ^m$, and any fixed $w\in M$, the assignment $u\mapsto f(w,u)$ is a diffeomorphism from $S$ onto the unit  sphere with respect to the {\em Riemannian} norm $\Vert ~\Vert _{w}$. So, for any nongeodesic 
$x\in F$ with $\Vert \nabla _tx^{(1)}(t)\Vert _{x(t)}$ constant, there are unique curves  
$y:[t_0,t_1]\rightarrow \R ^m$, $u:[t_0,t_1]\rightarrow S$ and $z:[t_0,t_1]\rightarrow (0,\infty )$ satisfying 
\begin{eqnarray}
\label{xeq}x^{(1)}(t)&=&y(t)\\
\label{yeq}y^{(1)}(t)&=&zf(x(t),u(t))  -\Gamma (x(t))(y(t),y(t))\\
\label{zeq}z^{(1)}(t)&=&0
\end{eqnarray}
\begin{equation}\label{bcons}\hbox{with }\quad 
x(t_0)~=~x_0,~~
x(t_1)=x_1,\quad \hbox{and}\quad 
y(t_0)~=~v_0,~~
y(t_1)~=~v_1.
\end{equation}
So an optimal curve $x$ corresponds to a {\em control extremal}, namely an extremal $(x,y,z):[t_0,t_1]\rightarrow M\times \R ^m\times (0,\infty )$ with control $u:[t_0,t_1]\rightarrow S$ for the optimal control problem with dynamics (\ref{xeq}), (\ref{yeq}), (\ref{zeq}), boundary conditions 
(\ref{bcons}), minimising 
$$(t_1-t_0)z~=~\int _{t_0}^{t_1}z(t)~dt.$$
The admissible controls are taken as the piecewise-continuous functions $u:[t_0,t_1]\rightarrow S$. This is equivalent to varying $x:[t_0,t_1]$ through piecewise-$C^2$ curves with $\Vert \nabla _tx^{(1)}\Vert _{x(t)}$ 
constant, while satisfying $x(t_i)=x_i$ and $x^{(1)}(t_i)=v_i$ for $i=0,1$.  

\spp
As a control extremal, $(x,y,z)$ and $u$ satisfy the Pontryagin Maximum Principle \cite{ponty}, stated in terms of the Hamiltonian  
$H:\{ 0,1\}\times (M\times \R ^m\times (0,\infty ) )\times (\R ^m\times \R ^m\times \R )^*\times S \rightarrow \R$ 
where   
$$H(\epsilon ,x,y,z,\lambda ,\mu,\nu,u)~:=~-\epsilon z+\lambda (y)  
+\mu (zf(x,u)-\Gamma _x(y,y)).$$
The Pontryagin Principle concerns a number $\epsilon \in \{ 0, 1\}$ and a curve of costates $(\lambda ,\mu ,\nu ):[t_0,t_1]\rightarrow (\R ^m\times \R ^m\times \R)^*$ associated with the control extremal $(x,y,z)$, satisfying 
\begin{eqnarray}
\label{lameq}\lambda ^{(1)}&=&\mu (z\langle \Gamma (e_j,\hat u),\hat u\rangle _x\hat u+d\Gamma _x(e_j)(y,y))e_j^*\\
\label{mueq}\mu^{(1)}&=&-\lambda +2\mu (\Gamma _x(e_j,y)) e_j^*\\
\label{nueq}\nu ^{(1)}&=&\epsilon -\mu (\hat u)
\end{eqnarray}
for almost all $t\in [t_0,t_1]$. Here we use Lemma \ref{lem3} to differentiate $\hat u:=f(x,u)$, $e_j^*$ is the dual of $e_j$, and we sum over $j=1,2,\ldots ,m$. Moreover, if $\epsilon =0$ then $(\lambda ,\mu ,\nu )$ is nowhere trivial. 

\spp
Corresponding to the absence of conditions on $z(t_0)$ and $z(t_1)$ are  transversality conditions  
\begin{equation}\label{trans}\nu (t_0)~=~\nu (t_1)~=~0.~\end{equation}
The Pontryagin Principle asserts that, for almost every $t\in [t_0,t_1]$,
 $u(t)$ is a maximiser of 
$H(\epsilon ,x(t),y(t),z,\lambda (t),\mu (t),\nu (t),u)$, namely whenever $\mu (t)\not= {\bf 0}$ 
$$u(t)~=~\frac{\tilde \mu (t)~~}{\Vert \tilde \mu (t)\Vert _E}.$$
\section{Local Proof of Theorem \ref{mainthm}}\label{localcalcsec}
We continue to suppose that $x$ is an optimal curve in an open subset $M$ of $\R ^m$.
\begin{lemma}\label{xlem} $\epsilon =1$ and $\Vert \tilde \mu (t)\Vert _{x(t)}$ has mean $1$. 
\end{lemma}

\spp
{\bf Proof of Lemma \ref{xlem}} By (\ref{nueq}) and transversality, 
$\int _{t_0}^{t_1}\Vert \tilde \mu (t)\Vert _{x(t)}~dt~=~\epsilon (t_1-t_0)$. 
So if $\epsilon =0$ we must have $\mu $ identically ${\bf 0}$. Then $\nu \equiv {\bf 0}$ by (\ref{nueq}) and transversality, and $\lambda \equiv {\bf 0}$ by (\ref{mueq}). But for $\epsilon =0$,
 $(\lambda ,\mu ,\nu )$ is nowhere trivial. So $\epsilon=1$ after all, and $\int _{t_0}^{t_1}\Vert \tilde \mu (t)\Vert _{x(t)}~dt =t_1-t_0$. So Lemma \ref{xlem} is proved.

\spp
\begin{lemma}\label{ylem} $L(\tilde \mu )$ is identically ${\bf 0}$.
\end{lemma}

\spp
{\bf Proof of Lemma \ref{ylem}:} For $t\in [t_0,t_1]$, rewrite (\ref{lameq}) and (\ref{mueq})  as 
\begin{eqnarray}
\label{newlameq}\nabla _t\lambda &=&(\mu (z\langle \Gamma (e_j,\hat u),\hat u\rangle _x\hat u+d\Gamma _x(e_j)(y,y)) -\lambda (\nabla _te_j ))e_j^*\\
\label{newmueq}\nabla _t\mu&=&-\lambda +\mu (\nabla _te_j) e_j^*.
\end{eqnarray}
Substituting for $\lambda $ from (\ref{newmueq}) in (\ref{newlameq}), 
\begin{equation}
\label{newnewlameq}\nabla _t\lambda =(\mu (z\langle \Gamma (e_j,\hat u),\hat u\rangle _x\hat u+d\Gamma _x(e_j)(y,y)- e_k^*(\nabla _te_j) \nabla _te_k) +\nabla _t(\mu )(\nabla _te_j)e_j^*.
\end{equation}
Differentiating (\ref{newmueq}) and substituting for $\nabla _t\lambda $ from (\ref{newnewlameq}), we find 
$$\nabla _t^2\mu ~=~-\nabla _t\lambda +(\nabla _t\mu )(\nabla _te_j)e_j^*+\mu (\nabla _t^2e_j) e_j^*
+\mu (\nabla _te_j) \nabla _te_j^*~=~$$
$$-\mu (z\langle \Gamma (e_j,\hat u),\hat u\rangle _x\hat u+d\Gamma _x(e_j)(y,y)-\nabla _t^2e_j) e_j^*.$$
Now ~$\displaystyle{d\Gamma _x(e_j)(y,y)
-\nabla _t^2e_j=\frac{\partial \Gamma ^k_{p,q}}{\partial x_j}y_{p}y_{q}e_k-\nabla _t(\Gamma ^k_{pj}y_{p}e_k)=}$ 
$$\frac{\partial \Gamma ^k_{p,q}}{\partial x_j}y_{p}y_{q}e_k-\frac{\partial \Gamma ^k_{pj}}{\partial x_q}y_{p}y_{q}e_k-\Gamma ^k_{pj}y_{p}^{(1)}e_k-
\Gamma ^r_{pj}\Gamma ^k_{qr}y_{p}y_{q}e_k~=~$$
$$\frac{\partial \Gamma ^k_{p,q}}{\partial x_j}y_{p}y_{q}e_k-(\frac{\partial \Gamma ^k_{pj}}{\partial x_q}-\Gamma ^k_{rj}\Gamma ^r_{q,p}+
\Gamma ^r_{pj}\Gamma ^k_{qr})y_{p}y_{q}e_k-z\Gamma ^k_{pj}\hat u_pe_k~=~R^k_{pjq}y_{p}y_{q}e_k-z\Gamma ^k_{pj}\hat u_pe_k$$
by (\ref{yeq}), (\ref{curveq}). Because $\mu (\hat u ) \hat u=\Vert \tilde \mu \Vert _x\hat u=\tilde \mu $, we have 
$\mu (\langle \Gamma (e_j,\hat u),\hat u\rangle _x\hat u-\Gamma ^k_{pj}\hat u_pe_k)=$
$$\langle \Gamma (e_j,\hat u),\hat u\rangle _x\Vert _x\tilde \mu \Vert _x-\mu (\Gamma (e_j,\hat u))~=~
\langle \Gamma (e_j,\hat u),\tilde \mu \rangle _x\Vert _x-\mu (\Gamma (e_j,\hat u))~=~0.$$
So ~$\displaystyle{\nabla _t^2\mu = 
-\mu (R^k_{pjq}y_{p}y_{q}e_k) e_j^*}$, ~and  
$$\displaystyle{\nabla _t^2\tilde \mu ~=~ 
-\langle \tilde \mu , R^k_{pjq}y_{p}y_{q}e_k\rangle _x e_j=-\langle \tilde R(e_j,y)y,\tilde \mu \rangle _xe_j}~=~-\langle R(\tilde \mu ,y)y,e_j \rangle _xe_j ~=~-R(\tilde \mu ,y)y $$
using symmetries of Riemannian curvature \cite{hawking} (2.27). This proves Lemma \ref{ylem}. 

\spp
So $\tilde \mu$ is a piecewise-$C^2$ field defined along the optimal curve $x$. Defining $\Phi :[t_0,t_1]\rightarrow [0,\infty )$ by 
$$\Phi (t)~:=~\Vert \tilde \mu (t)\Vert _{x(t)}^2,$$
we have $L(\sqrt{\Phi (t)}\nabla _tx^{(1)})=L(z\tilde \mu )=zL(\tilde \mu )={\bf 0}$, since $z$ is constant and by Lemma \ref{ylem}.   
\begin{lemma}\label{lem4} $\Phi ^{-1}(0)$ is a discrete subset of $[t_0,t_1]$.
\end{lemma}

\spp
{\bf Proof of Lemma \ref{lem4}:} Let  $s_*$ be an accumulation point of $\Phi ^{-1}(0)$. Then $\Phi (s_*)=\Phi ^{(1)}(s_*)=\Phi ^{(2)}(s_*)=0$ (the second derivative is taken to be one-sided if $s_*$ is  one of the finitely many points of discontinuity of $\Phi ^{(2)}$). So  
$$\frac{d^2}{dt^2}\langle \tilde \mu ,\tilde \mu \rangle _{x(s_*)}~=~2\frac{d}{dt}\langle \nabla _t\tilde \mu ,\tilde \mu \rangle _{x(s_*)}~=~2\langle \nabla _t^2\tilde \mu ,\tilde \mu \rangle _{x(s_*)}+2\langle \nabla _t\tilde \mu ,\nabla _t\tilde \mu \rangle _{x(s_*)}~=~0.$$
So $\tilde \mu (s_*)=(\nabla_t \tilde \mu )(s_*)={\bf 0}$. Since $L$ is a second order linear differential operator, $\tilde \mu $ is identically ${\bf 0}$, by Lemma \ref{ylem}. This contradiction of Lemma \ref{xlem} proves Lemma \ref{lem4}.

\spp
Since $\tilde \mu $ is piecewise-$C^2$, so is $\nabla _tx^{(1)}$ on any open interval $I\subset [t_0,t_1]$ where $\Phi $ is nonzero. So $x\vert I$ is piecewise-$C^4$. So $\tilde \mu \vert I$ is also piecewise-$C^4$, and so on. In fact $x\vert I$ and $\tilde \mu \vert I$ are $C^\infty$. So $\Phi ^{(2)}$ is continuous except where $\tilde \mu$ is ${\bf 0}$.  So $x$ and $\tilde \mu $ are $C^\infty$ except where $\Phi $ is $0$. 

\spp
Theorem \ref{mainthm} is now proved for the special case where $M$ is open in $\R ^m$. As argued at the end of \S \ref{Linfsec},  this completes the proof of Theorem \ref{mainthm} in general. 
\section{Spheres}\label{spheresec}
Take $M$ to be the unit $m$-sphere $S^m$ in $E^{m+1}$. The Levi-Civita covariant derivative $\nabla _YZ$ of a vector field $Z$ on $S^m$ in the direction of $Y\in TS^m_y=y^\perp$ is given by 
$$\nabla _YZ~=~dZ_y(Y)+\langle Y,Z\rangle y$$
where $\langle ~,~\rangle $ is the Euclidean inner product, $Z$ is treated as a map from $S^m$ to $E^{m+1}$, and $y\in S^m$. Accordingly, the Riemannian curvature $R$ on $S^m$ is given by 
$$R(X,Y)Z~=~\langle Y,Z\rangle X-\langle X,Z\rangle Y\quad \hbox{and, for a field }X~ \hbox{along} ~x:[t_0,t_1]\rightarrow S^m,$$ 
$$\displaystyle{L(X)~=~X^{(2)}+\langle x^{(1)},x^{(1)}\rangle X+(\langle x^{(2)},X\rangle +2\langle x^{(1)},X^{(1)}\rangle )x}.$$
If $x$ is optimal then, by Theorem \ref{mainthm}, for some $z\in \R$ and some field $X$ with a discrete set of zeroes, we have 
\begin{eqnarray}
\label{xeq1}x^{(2)}&=&zX/\Vert X\Vert -\langle x^{(1)},x^{(1)}\rangle x\\
\label{Xeq}X^{(2)}&=&-\langle x^{(1)},x^{(1)}\rangle X-(\langle x^{(2)},X\rangle +2\langle x^{(1)},X^{(1)}\rangle )x
\end{eqnarray}
wherever $X(t)\not= {\bf 0}$.
\begin{example}\label{sphereex} Take $m=2$, $z=1.2$, and 
$$x(0)=(1,0,0),~x^{(1)}(0)=(0,1,0),~ 
X(0)=(0,1,200),~X^{(1)}(0)=(-1,2,1).$$
The solution $x:[0,8]\rightarrow S^2$ is shown, viewed from $(2,2,-1)$, in Figure \ref{xfig}, beginning at the left, and ending at $x(8)\approx (-0.433207, 0.898726, 0.0679917)$.  
\begin{figure}[htbp] 
   \centering
   \includegraphics[width=2in]{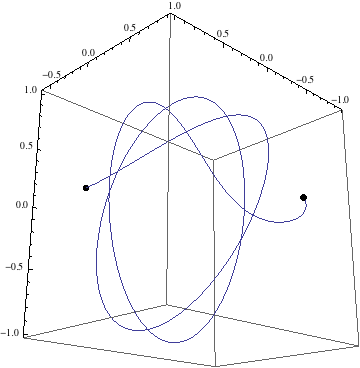} 
   \caption{$x:[0,8]\rightarrow S^2$ in Example \ref{sphereex}}
   \label{xfig}
\end{figure}
\end{example}
\section{Bi-Invariant Lie Groups}\label{liegp}
Take $M$ to be a Lie group with  a bi-invariant Riemannian metric $\langle ~,~\rangle $, as in \cite{milnor} \S 21. The {\em (left) Lie-reduction} of a field $Y$ along $x:[t_0,t_1]\rightarrow G$ is the curve $Y_L$ in the Lie algebra ${\cal G}$ of $G$ given by 
$$Y_L(t)~:=~dL_{x(t)}^{-1}Y_{x(t)}$$
where $L(g):G\rightarrow G$ denotes left-multiplication by $g\in G$. 
Given $x_0,x_1,v_0,v_1$, and an $F_{x_0,x_1,v_0,v_1,t_0,t_1}$-optimal curve $x:[t_0,t_1]\rightarrow G$, denote the left Lie reduction of $x^{(1)}$ by $V$. By Theorem \ref{mainthm}, $z=\Vert V^{(1)}(t)\Vert $ is constant, and 
$X_L=\phi V^{(1)}$ where $X:=\phi \nabla _tx^{(1)}$. 

\spp
Since $X_L$ is a multiple of $V^{(1)}$, comparison of (\ref{eleq}) with Lemmas 1, 2 of \cite{jmp} gives
\begin{equation}\label{reducedeleq}
\tilde X_L^{(2)}~=~[\tilde X_L^{(1)},V]
\end{equation}
namely $X_L^{(1)}$ defines a Lax constraint \cite{lyledual} on $V$. So by Theorem 2.1 of \cite{lyleqjm}, if $G$ is semisimple, $x$ is generically obtainable by quadrature from $X_L^{(1)}=\phi ^{(1)}V^{(1)}+\phi V^{(2)}$ and $V$. Also 
$c:=\Vert X_L^{(1)}\Vert ^2=z^2\phi ^{(1)}(t)^2+\phi (t)^2\Vert V^{(2)}(t)\Vert ^2$
is constant, by (\ref{reducedeleq}). Integrating (\ref{reducedeleq}),
\begin{equation}\label{intreducedeleq}
X_L^{(1)}~=~[X_L,V]+C
\end{equation} 
where $C\in {\cal G}$ is constant. Taking inner products of both sides of (\ref{intreducedeleq}) with $V^{(1)}$,
$$z\phi ^{(1)}=\langle C,V^{(1)}\rangle ~\Longrightarrow ~z\phi (t)=\langle C,V(t)\rangle +a$$
where $a$ is constant. 
Taking inner products of both sides of (\ref{intreducedeleq}) with $V$,
$$\langle X_L^{(1)},V\rangle ~=~\langle C,V\rangle ~=~z\phi -a .$$

\spp
By analogy with Riemannian cubics \cite{jmp}, the optimal curve $x$ is called 
{\em null} when $C={\bf 0}$. 
Then, since $X_L$ defines another Lax constraint on $V$, $\Vert X_L\Vert ^2=z^2\phi (t)^2$ is constant. So $\phi $ is constant unless $x$ is geodesic. Consequently,  by equation (\ref{eleq}) 
in Theorem \ref{mainthm}, we have the 
\begin{corollary} For $M$ a bi-invariant Lie group, $x\in F_{x_0,x_1,v_0,v_1,t_0,t_1}$ is a null optimal curve if and only if $x$ is a null Riemannian cubic. 
\end{corollary}
Null Riemannian cubics are studied in \cite{jmp}, \cite{lyledual}, \cite{lyleSIAM} for $G=SO(3)$, 
but most optimal curves are non-null and have different geometry, as seen from their reductions. 
\begin{example}\label{so3ex} Take $G=SO(3)$, with the bi-invariant Riemannian metric for which 
$${\rm ad}:E^3\rightarrow so(3)$$
is an isometry, where Euclidean $3$-space $E^3$ is viewed as a Lie algebra with respect to the cross-product $\times $. Then the Lie reduction of an optimal curve $x^{(1)}$ can be identified with a curve $V$ in $E^3$, satisfying 
\begin{eqnarray}
\label{veq}V^{(1)}(t)&=&zW/\Vert W\Vert \\
\label{weq}W^{(1)}(t)&=&W\times V+C 
\end{eqnarray}
where $z\in \R$, $C\in E^3$, $\Vert ~\Vert $ is the Euclidean norm $\Vert ~\Vert _E$, and $W(t):=X_L(t)\in E^3$. As noted previously, there are conserved quantities
$$c~=~\Vert W^{(1)}\Vert ^2,\quad a~=~z\Vert W\Vert -\langle C,V\rangle ~=~z\Vert W\Vert -\langle W^{(1)},V\rangle .$$
Taking $z=1.2$, $C=-(2,1,0)$, $V(0)=(1,2,3)$, $W(0)=(-1,-4,6)$, a numerical solution of (\ref{veq}), (\ref{weq}) for $V:[0,700]\rightarrow E^3$ is shown in Figure \ref{vfig} below. 
The inner dot is at $V(0)$ and the outer is at the other endpoint $V(700)\approx 
(2.36765, 4.69752, 8.40276)$, and the image of $V$ is viewed from $(10,8,0)$. The curve is nonplanar and seems to be approximately periodic. 
\begin{figure}[htbp] 
   \centering
   \includegraphics[width=2in]{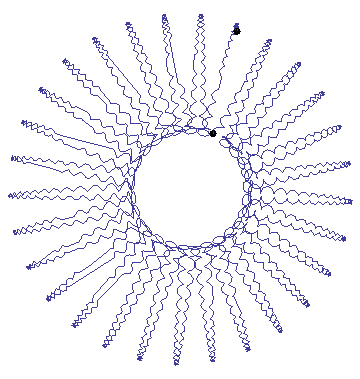} 
   \caption{Lie Reduction of $x^{(1)}:[0,700]\rightarrow SO(3)$ in Example \ref{so3ex}}
   \label{vfig}
\end{figure}
Similarly, $\phi =\Vert W\Vert :[0,700]\rightarrow \R$  seems to be approximately periodic, and is apparently everywhere nonzero. The graph over the interval $[0,55]$ is shown in Figure \ref{phfig}. 
\begin{figure}[htbp] 
   \centering
   \includegraphics[width=3in]{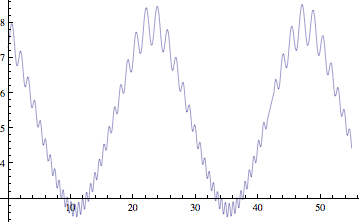} 
   \caption{$\phi $ in Example \ref{so3ex}}
   \label{phfig}
\end{figure}
If $C$ is replaced by $(2,1,0)$, the reduction $V$ of $x^{(1)}$ appears to be unbounded. The image of $V\vert [0,5]$ shown in Figure \ref{vfig2} resembles that of a null Lie quadratic, spiralling curve inwards towards an asymptotic line. However $x$ is not a null Riemannian cubic, because $C\not={\bf 0}$. The lower dot corresponds to $V(0)$ and the upper dot is at $V(5)\approx (1.77133, 4.50895, 7.05963)$. 
\begin{figure}[htbp] 
   \centering
   \includegraphics[width=2in]{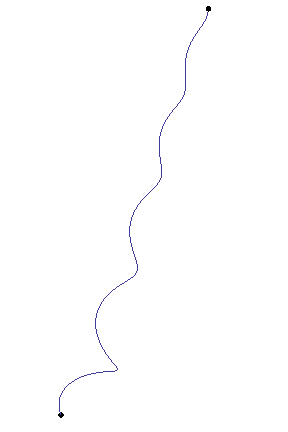} 
   \caption{A nonperiodic reduction of a non-null optimal curve in Example \ref{so3ex}}
   \label{vfig2}
\end{figure} 
\end{example}
\section{Two Points and One Velocity}\label{2pts1velsec}
So far, for feasible curves $x:[t_0,t_1]\rightarrow M$, $x(t_i)$ and $x^{(1)}(t_i)$ have been prescribed for $i=0,1$. With no conditions on the $x^{(1)}(t_i)$, $J$ would be minimised by a geodesic. Consider next what happens when the condition on $x^{(1)}(t_1)$ is lifted. 
\begin{definition} Let $x_0,x_1\in M$ and $v_0\in TM_{x_0}$ be given. Let $F_{x_0,x_1,v_0,t_0,t_1}$ be the set of piecewise-$C^2$ curves for which $x(t_i)=x_i$ for $i=0,1$ and $x^{(1)}(t_0)=v_0$. Defining 
$J:F_{x_0,x_1,v_0,t_0,t_1}\rightarrow \R$ in the same way as before, a minimiser $x\in F_{x_0,x_1,v_0,t_0,t_1}$ of $J$ is said to be {\em $F_{x_0,x_1,v_0,t_0,t_1}$-optimal}.  
\end{definition}
\begin{theorem}\label{subthm1} Let $x$ be $F_{x_0,x_1,v_0,t_0,t_1}$-optimal. Then $\Vert  \nabla _tx^{(1)} \Vert _{x(t)}$ is constant. For some piecewise-$C^2$ function 
$\Phi :[t_0,t_1]\rightarrow [0,\infty )$ with $\Phi ^{-1}(0)$ a discrete set containing $t_1$, we have 
\begin{equation}\label{eleq1}L(\phi \nabla _tx^{(1)})~=~{\bf 0}\end{equation}
where $\phi (t):=\sqrt{\Phi (t)}$. Moreover $x$ and $\phi $ are $C^\infty$ except possibly where $\Phi (t)=0$. 
\end{theorem}

\spp
The change in the conclusion from Theorem \ref{mainthm} is that $\Phi (t_1)=0$. We omit the changes of detail needed for the proof of Theorem \ref{subthm1}. Considering track-inverses, a similar result holds when $x_0,x_1\in M$ and $v_1\in TM_{x_1}$ are given, and feasibility means $x(t_i)=x_i$ for $i=0,1$ with $x^{(1)}(t_1)=v_1$. 
\section{Many Points and No Velocities}\label{manyptssec}
Another kind of feasibility occurs when, for $n\geq 3$, we are given real numbers $t_0<t_1<\ldots <t_j<\ldots <t_n$ and points $x_0,x_1,\ldots ,x_j,\ldots ,x_n\in M$. Denote by $F_{x_0,x_1,\ldots ,x_n,t_0,t_1,\ldots ,t_n}$ the set of piecewise-$C^2$ curves 
$x:[t_0,t_n]\rightarrow M$  satisfying $x(t_i)=x_i$ for $i=0,1,\ldots ,n$. Define $J:F_{x_0,x_1,\ldots ,x_n,t_0,t_1,\ldots ,t_n}\rightarrow \R$ by 
$$J(x)~:=~\max \{  \Vert \nabla _tx^{(1)}\Vert _{x(t)}:t_0\leq t\leq t_n  \} .$$
Then a minimiser $x\in F_{x_0,x_1,\ldots ,x_n,t_0,t_1,\ldots ,t_n}$ of $J$ is called 
{\em $F_{x_0,x_1,\ldots ,x_n,t_0,t_1,\ldots ,t_n}$-optimal}. On considering optimality of track-summands, we find from Theorems \ref{mainthm}, \ref{subthm1},     
\begin{corollary}\label{subthm2} Let $x$ be 
$F_{x_0,x_1,\ldots ,x_n,t_0,t_1,\ldots ,t_n}$-optimal. Then, for some $j=1,2,\ldots ,n$ and all $t\in [t_{j-1},t_j]$, $\Vert \nabla _tx^{(1)}\Vert _{x(t)}$ is constant. Moreover, 
\begin{itemize}
\item for some piecewise-$C^2$  function $\Phi :[t_{j-1},t_j]\rightarrow [0,\infty)$ with $\Phi ^{-1}(0)$ discrete, we have 
$L(\phi \nabla _tx^{(1)})={\bf 0}$ for $t\in [t_{j-1},t_j]$, where $\phi (t):=\sqrt {\Phi (t)}$,
\item $x$ and $\phi $ are $C^\infty$ except possibly where $\Phi (t)=0$,
\item if $j=1$ then $\Phi (t_0)=0$,
\item if $j=n$ then $\Phi (t_n)=0$.
\end{itemize}  
\end{corollary}

\spp
{\bf Proof of Corollary:} Suppose there is no such $j$. Then, by Theorems \ref{mainthm}, \ref{subthm1}, for every $j=1,2,\ldots ,n$, $x\vert [t_{j-1},t_j]$ is suboptimal with respect to the restriction of $J_\infty$ to $F_{x(t_{j-1}),x(t_j),x^{(1)}(t_{j-1}),x^{(1)}(t_j)}$. This contradicts optimality of $x$.

\spp
\begin{example} For a complete Riemannian 
manifold $M$, 
 choose $x_0,x_1,x_2\in M$ not lying on a geodesic. Choose real numbers $t_0<t_1<t_2$ and let  $x$ be $F_{x_0,x_1,x_2,t_0,t_1,t_2}$-optimal.  Set $\varepsilon :=J_\infty (x)$. Then $\varepsilon >0$. 
 
 \spp
For some $t_3>t_2$ let $y:[t_2,t_3]\rightarrow M$ be a geodesic and set $x_3:=y(t_3)$. Let $\tilde y:[t_2,t_3]\rightarrow M$ be be any piecewise-$C^2$ curve with $\tilde y^{(1)}(t_2)=x^{(1)}(t_2)$, $\tilde y(t_3)=x_3$ and $\Vert \nabla _t\tilde y^{(1)}\Vert _{\tilde y(t)}<\varepsilon $ for all $t\in (t_2,t_3)$. For instance $\tilde y$ could be chosen as a perturbation of $y$, or indeed $y$ itself.

\spp
Then the track-sum $\bar x:[t_0,t_3]\rightarrow M$ of $x$ and $\tilde y$ lies in $ F_{x_0,x_1,x_2,x_3,t_0,t_1,t_2,t_3}$. Moreover $J_\infty (\bar x)=J_\infty (x)$, and $\bar x$ is easily seen to be $F_{x_0,x_1,x_2,x_3,t_0,t_1,t_2,t_3}$-optimal. Yet, since there are continuously many 
possibilities for $\bar x$, the conditions of Corollary \ref{subthm2} do not apply for $j=3$.  
\end{example}
\section{Conclusion}
This paper proves necessary conditions for minimum $L^\infty$ acceleration curves in Riemannian manifolds $M$, analogous to the conditions in \cite{lylegreg} for Riemannian cubics. 
Results of this kind were previously known only when the manifold $M$ was flat. In the present paper the necessary conditions are examined in detail when $M$ is a sphere or a bi-invariant Lie group. Examples are given in the case where $M$ is $S^2$ or $SO(3)$, raising  questions about the asymptotics and symmetry of optimal curves. 
\section{Acknowledgements}
I am very grateful to Professor Charles Micchelli for bringing \cite{charlie} to my attention, for suggesting to use the $L^\infty$ extremals of \S \ref{Linfsec} in a Riemannian context, and for his generous hospitality. 
The use of optimal control is also essential in the proof of Theorem \ref{mainthm}. In continuing research, Professor Yalcin Kaya and I applied the Pontryagin Principle in a different way  to related problems in a Euclidean context. I am very grateful to Professor Kaya for stimulating conversations and for his kind hospitality. 
%





\bibliographystyle{model1-num-names}
\bibliography{<your-bib-database>}

\begin{thebibliography}{00}
\bibitem{agrachev}
A.A. Agrachev and Y.L. Sachkov, 
{\em Control Theory From The Geometric Viewpoint}, 
Encyclopaedia of Math. Sciences, Control Theory and Optimization {\bf 11}, Springer (2004). 

\bibitem{barbero}
M. Barbero-Linan and M.C. Munoz-Lecanda,
``Geometric approach to Pontryagin's maximum principle,''
{\em Acta Appl. Math.} {\bf 108} (2009) 429--485.

\bibitem{bloch}
A.M. Bloch with J. Baillieul, P. Crouch and J. Marsden, 
{\em Nonholonomic Mechanics and Control}, 
Interdisciplinary Applied Mathmatics, Springer (2003).

\bibitem{chang}
D.E. Chang,
``A simple proof of the Pontryagin maximum principle on manifolds,''
{\em Automatica} {\bf 47} (2011) 630--633.

\bibitem{clarke}
F.H. Clarke and R.B. Vinter,
``Optimal multiprocesses,''
{\em SIAM J. Control \& Optimization} {\bf 27} (5)  (1989) 1072--1090.

\bibitem{crouch1} 
M. Camarinha, F. Silva Leite, P. Crouch,
``On the geometry of Riemannian cubic polynomials,'' 
{\em Differential Geom. Appl.} {\bf 15} (2) (2001) 107--135.

\bibitem{crouch4}
P. Crouch and F. Silva Leite,  
``The dynamic interpolation problem: on Riemannian manifolds, Lie groups, 
and symmetric spaces,''
{\em J. Dynam. Control Systems} {\bf 1} (2) (1995) 177--202. 


\bibitem{deboorbook}
Carl de Boor,
{\em A Practical Guide to Splines}, 
Applied Mathematical Sciences {\bf 27}, Springer (2001).

\bibitem{deboor}
Carl de Boor,
``A remark concerning perfect splines,''
{\em Bull. Amer. Math. Soc.} {\bf 80} (4) (1974), 724--727.

\bibitem{jerome}
S.D. Fisher and J.W. Jerome,
``The existence, characterization and essential uniqueness of solutions of $L^\infty$ extremal problems,''
{\em Trans. Amer. Math. Soc.} {\bf 187} (1974) 391--404.

\bibitem{glaeser}
G. Glaeser,
``Prolongement extremal de fonctions differentiables,''
{\em Publ. Sect. Math. Facult\'e des Sciences Rennes}, Rennes, France (1967).

\bibitem{giambo}
R. Giambo, F. Giannoni and P. Piccione, 
``An analytical theory for Riemannian cubic polynomials," 
{\em IMA J. Math. Control \& Information} {\bf 19} (2002) 445--460.

\bibitem{hawking}
S.W. Hawking and G.F.R. Ellis,
{\em The Large Scale Structure of Space-Time},
Cambrdge Monographs on Mathematical Physics,
Cambridge UP (1973). 


\bibitem{vega}
S. Guti\'errez, J. Rivas and L. Vega,
``Formation of singularities and self-similar vortex motion under the localized induction 
approximation,"
{\em Comm. in Partial Differential Equations} {\bf 28} (2003) 927--968. 



\bibitem{charlie}
C.A. Micchelli,
``Curves from variational principles,''
{\em RAIRO Mod\'elisation Math\'ematique et Analyse Num\'erique} {\bf 26} (1) (1992) 77--93.

\bibitem{milnor}
J. Milnor,
``Morse Theory,''
{\em Annals of Math. Studies} {\bf 51}, Princeton UP (1963). 

\bibitem{lylegreg}
L. Noakes, G. Heinzinger and B. Paden,
``Cubic splines on curved spaces,''
{\em IMA J. Math. Control \& Information} {\bf 6} (1989) 465--473.

\bibitem{jmp}
L. Noakes,
``Null cubics and Lie quadratics,''
{\em J. Math. Physics} {\bf 44} (3) (2003) 1436--1448. 

\bibitem{lylenonnull}
L. Noakes,
``Non-null Lie quadratics in $E^3$,'' 
{\em J. Math Physics}, {\bf 45} (11) (2004) 4334--4351.

\bibitem{lyledual}
L. Noakes,
``Duality and Riemannian Cubics,'' 
{\em Adv. in Computational Math.} {\bf 25} (2006) 195--209.

\bibitem{lyleqjm}
L. Noakes,
``Lax constraints in semisimple Lie groups,''
{\em Quart. J. Math.}, {\bf 57} (2006) 527--538. 

\bibitem{lyleSIAM}
L. Noakes, 
``Asymptotics of null Lie quadratics in $E^3$,'' 
{\em SIAM J. on Applied Dynamical Systems}, {\bf 7} (2) (2008) 437--460. 

\bibitem{lyletomaszrev}
L. Noakes and T. Popiel,
``Geometry for robot path planning," 
{\em Robotica} {\bf 25} (2007) 691--701.

\bibitem{pauley}
M. Pauley, 
``Cubics, Curvature and Asymptotics,''
{\em PhD Thesis, University of Western Australia}, (2011).   


\bibitem{ponty}
L.S. Pontryagin, V.G. Boltyanski, R.V. Gamkrelidze and E.F. Mischenko, 
(translated by K.N. Trirogoff),
{\em The Mathematical Theory of Optimal Processes},
Interscience, John Wiley (1962).

\bibitem{popiel}
T. Popiel,
``Geometrically-Defined Curves in Riemannian Manifolds,''
{\em PhD Thesis, University of Western Australia}, (2007).



\bibitem{singer}  
D.A. Singer,
``Lectures on elastic curves and rods,'' 
{\em Curvature and Variational Modeling in Physics and Biophysics,  AIP Conf. Proc.} {\bf 1002} (2008) 3--32. 
%



\end{thebibliography}



\end{document}